\documentclass[11pt]{article}
\usepackage{booktabs}
\usepackage{caption}
\usepackage{mathrsfs}
\usepackage{amsmath}
\usepackage{amsfonts,amsthm,amssymb,mathrsfs,bbding}
\usepackage{float}
\usepackage{graphics,epstopdf}
\usepackage{graphics,multicol}
\usepackage{graphicx}
\usepackage{color}
\usepackage{enumerate}
\usepackage{tikz}
\usepackage{caption}
\usepackage{rotating}
\usepackage{lscape}
\usepackage{longtable}
\allowdisplaybreaks[4]
\usepackage{tabularx}
\usepackage[colorlinks=true,anchorcolor=blue,filecolor=blue,linkcolor=blue,urlcolor=blue,citecolor=blue]{hyperref}
\usepackage{extarrows}
\usepackage{cite}
\usepackage{latexsym,bm}
\usepackage{mathtools}
\evensidemargin=\oddsidemargin\topmargin=-1.5cm

\newtheorem{thm}{Theorem}[section]
\newtheorem{prob}[thm]{Problem}

\newtheorem{lem}{Lemma}[section]
\newtheorem{cor::3.1}{Corollary}[section]

\newtheorem{remark}[thm]{Remark}

\newtheorem{conj}{Conjecture}

\theoremstyle{definition}
\renewcommand\proofname{\bf Proof}
\addtocounter{section}{0}
\begin{document}

\title{\LARGE{\bf Toughness, hamiltonicity and spectral radius in graphs} \footnote{This work is supported by the National Natural Science Foundation of China (Grant Nos. 11771141, 12011530064 and 11871391)}\setcounter{footnote}{-1}\footnote{\emph{Email address:} ddfan0526@163.com (D. Fan), huiqiulin@126.com (H. Lin), luhongliang@mail.xjtu.edu.cn (H. Lu).}}
\author{Dandan Fan$^{a,b}$,  Huiqiu Lin$^a$\thanks{Corresponding author.}, Hongliang Lu$^c$\\[2mm]
\small\it $^a$ School of Mathematics, East China University of Science and Technology, \\
\small\it   Shanghai 200237, China\\[1mm]
\small\it $^b$ College of Mathematics and Physics, Xinjiang Agricultural University\\
\small\it Urumqi, Xinjiang 830052, China\\[1mm]
\small\it $^c$ School of Mathematics and Statistics, Xi'an Jiaotong University\\
\small\it Xi'an, Shaanxi 710049, China}
\date{}
\maketitle

{\flushleft\large\bf Abstract}
The study of the existence of hamiltonian cycles in a graph is a classic problem in graph theory. By incorporating toughness and spectral conditions, we can consider Chv\'{a}tal's conjecture from another perspective: what is the spectral condition to guarantee the existence of a hamiltonian cycle among $t$-tough graphs?
We first give the answer to $1$-tough graphs, i.e. if $\rho(G)\geq\rho(M_{n})$, then $G$ contains a hamiltonian cycle, unless $G\cong M_{n}$, where $M_{n}=K_{1}\nabla K_{n-4}^{+3}$ and $K_{n-4}^{+3}$ is the graph obtained from $3K_{1}\cup K_{n-4}$ by adding three independent edges between $3K_{1}$ and $K_{n-4}$.
The Brouwer's toughness theorem states that every $d$-regular connected graph always has $t(G)>\frac{d}{\lambda}-1$ where $\lambda$ is the second largest absolute eigenvalue of the adjacency matrix.
In this paper, we extend the result in terms of its spectral radius, i.e. we provide a spectral condition for a graph to be 1-tough with minimum degree $\delta$ and to be $t$-tough, respectively.

\begin{flushleft}
\textbf{Keywords:} Toughness; hamiltonian cycle; spectral radius.
\end{flushleft}
\textbf{AMS Classification:} 05C50

\section{Introduction}

A graph $G$ is \textit{$t$-tough} if $|S|\geq tc(G-S)$ for every subset $S\subseteq V(G)$ with $c(G-S)>1$, where $c(G)$ is the number of components of a graph $G$.
The \textit{toughness $t(G)$} of $G$ is the maximum $t$ for which $G$ is $t$-tough (taking $t(K_{n})=\infty$). Hence if $G$ is not a complete graph, $t(G)=\min\{\frac{|S|}{c(G-S)}\}$ where the minimum is taken over all cut sets of vertices in $G$.
This concept was first introduced by
Chv\'{a}tal \cite{V.C} as a way of measuring how tightly various pieces of a graph hold together. 1-tough is a necessary condition for a graph to be hamiltonian.

Over the past decade, researchers are focused on finding spectral conditions for a graph to be hamiltonian. In 2010, Fiedler and Nikiforov \cite{M.F} first described several sufficient conditions, in terms of the eigenvalues of the adjacency matrices, for a graph to contain a hamiltonian cycle. Up to now, much attention has been paid on this topic, and we refer the reader to \cite{V.B,M.L,B.N}. By imposing the minimum degree of a graph as a parameter, Li and Ning \cite{B.L} extended the results in \cite{M.F,B.N}. Denote by $\delta(G)$ ($\delta$ for short) and $\Delta(G)$ the minimum degree and the maximum degree of $G$, respectively. Let $A(G)$ be the adjacency matrix of a graph $G$. The largest eigenvalue of $A(G)$, denoted by $\rho(G)$, is called the \textit{spectral radius} of $G$. Denote by $'\nabla'$ and $'\cup'$ the join and union product, respectively.
Note that $\delta\geq 2$ is a trivial necessary condition for a 1-tough graph. Inspired by the work of Li and Ning \cite{B.L}, it is interesting to find a spectral condition for a graph to be 1-tough for $\delta\geq 2$.
 \begin{thm}\label{thm::1.1}
Suppose that $G$ is a connected graph of order $n\geq \max\{5\delta, \frac{2}{5}\delta^{2}+\delta\}$ with $\delta\geq 2$. If $\rho(G)\geq\rho(K_{\delta}\nabla(K_{n-2\delta}\cup\delta K_{1}))$, then $G$ is a 1-tough graph, unless $G\cong K_{\delta}\nabla(K_{n-2\delta}\cup\delta K_{1})$.
\end{thm}

It is easy to see that if $G$ has a hamiltonian cycle then $G$ is 1-tough. Conversely, Chv\'{a}tal \cite{V.C} proposed the following conjecture.

\begin{conj}(See \cite{V.C})\label{conj::1.1}
There exists a finite constant $t_{0}$ such that every $t_0$-tough graph is hamiltonian.
\end{conj}

Around this conjecture, Enomoto, Jackson, Katerinis and Saito \cite{H.E-1} obtained
that every 2-tough graph contains a 2-factor and there exists $(2-\varepsilon)$-tough graphs without a 2-factor, where $\varepsilon>0$ is an arbitrary small number, and hence without a hamiltonian cycle. They conjectured that the value of $t_{0}$ might be 2. In 2000, this result was disproved by Bauer, Broersma and Veldman \cite{D.B-3}, and they observed that if such a $t_0$ exists, then it must be at least $\frac{9}{4}$. In general, the conjecture still remains open.

Another stream involved toughness is to find Dirac-type and Ore-type conditions for the existence of hamiltonian cycles in a graph. In 1995, Bauer, Broersma, van den Heuvel and Veldman \cite{DB-3} proved that a $t$-tough graph $G$ on $n\geq 3$ orders satisfying $\delta>\frac{n}{t+1}-1$ is hamiltonian.  Recently, Shan \cite{S.S} showed that a $t$-tough graph to be hamiltonian if the degree sum of any two non-adjacent vertices of $G$ is greater than $\frac{2n}{t+1}+t-2$. For more details, we refer the reader to \cite{D.B,D.B-1,B.F}. Naturally, we consider a graph to be hamiltonian by incorporating toughness and spectral conditions of the graph and raise the following problem.

\begin{prob}\label{prob1}
What is the spectral condition to guarantee the existence of a hamiltonian cycle among $t$-tough graphs?
\end{prob}
Based on the Dirac-type and Ore-type conditions for the 1-tough graphs, in this paper, we determine the spectral condition to guarantee the existence of a hamiltonian cycle among $1$-tough graphs. Let $K_{n-4}^{+3}$ be the graph obtained from $3K_{1}\cup K_{n-4}$ by adding three independent edges between $3K_{1}$ and $K_{n-4}$, and let $M_{n}=K_{1}\nabla K_{n-4}^{+3}$.

\begin{thm}\label{thm::1.2}
Suppose that $G$ is a connected 1-tough graph of order $n\geq 18$ with $\delta\geq 2$. If $\rho(G)\geq\rho(M_{n})$, then $G$ contains a hamiltonian cycle, unless $G\cong M_{n}$.
\end{thm}

Let $\lambda_i(G)$ denote the $i$th largest eigenvalue of $A(G)$.
By the Perron-Frobenius Theorem, $\lambda_1$ is always
positive (unless $G$ has no edges) and $|\lambda_i|\leq \lambda_1$ for all $i\geq 2$. Let $\lambda = \max_{2\leq i\leq n}|\lambda_i| = \max\{|\lambda_2|, |\lambda_n|\}$, where $n$ is the number of vertices of $G$. That is, $\lambda$ is the second largest absolute eigenvalue.

The study of toughness from eigenvalues was initiated by Alon \cite{Alon} who showed that for any connected $d$-regular graph $G$, $t(G)>\frac{1}{2}(\frac{d^2}{d\lambda+\lambda^2}-1)$. Around the same time, Brouwer \cite{Brouwer1} independently discovered a slightly better bound $t(G)>\frac{d}{\lambda}-2$, and he \cite{Brouwer1,Brouwer2} further conjectured that the lower bound can be improved to $t(G)>\frac{d}{\lambda}-1$. Subsequently, Gu \cite{Gu1} strengthened the result of Brouwer and showed the lower bound can be improved to $t(G)>\frac{d}{\lambda}-\sqrt{2}$. Very recently, Gu \cite{Gu2}  confirmed the conjecture of Brouwer. See more results on the relationship between the toughness and eigenvalues \cite{CG,CW,GH}. It is interesting to extend the results on finding a condition for a graph to be $t$-tough in terms of its spectral radius. In this paper, we give a spectral condition to guarantee the graph to be $t$-tough.

\begin{thm}\label{thm::1.3}
Let $t$ be a positive integer. If $G$ is a connected graph of order $n\geq 4t^2+6t+2$ with $\rho(G)\geq\rho(K_{2t-1}\nabla(K_{n-2t}\cup K_{1}))$, then $G$ is a $t$-tough graph, unless $G\cong K_{2t-1}\nabla(K_{n-2t}\cup K_{1})$.
\end{thm}

\section{Proof of Theorem \ref{thm::1.1}}
Let $M$ be a real $n\times n$ matrix, and let $X=\{1,2,\ldots,n\}$. Given a partition $\Pi=\{X_1,X_2,\ldots,X_k\}$ with $X=X_{1}\cup X_{2}\cup \cdots \cup X_{k}$, the matrix $M$ can be partitioned as
$$
M=\left(\begin{array}{ccccccc}
M_{1,1}&M_{1,2}&\cdots &M_{1,k}\\
M_{2,1}&M_{2,2}&\cdots &M_{2,k}\\
\vdots& \vdots& \ddots& \vdots\\
M_{k,1}&M_{k,2}&\cdots &M_{k,k}\\
\end{array}\right).
$$
The \textit{quotient matrix} of $M$ with respect to $\Pi$ is defined as the $k\times k$ matrix $B_\Pi=(b_{i,j})_{i,j=1}^k$ where $b_{i,j}$ is the  average value of all row sums of $M_{i,j}$.
The partition $\Pi$ is called \textit{equitable} if each block $M_{i,j}$ of $M$ has constant row sum $b_{i,j}$.
Also, we say that the quotient matrix $B_\Pi$ is \textit{equitable} if $\Pi$ is an equitable partition of $M$.

\begin{lem}(Brouwer and Haemers \cite[p. 30]{BH})\label{lem::2.1}
Let $M$ be a real symmetric matrix, and let $\lambda_{1}(M)$ be the largest eigenvalue of $M$. If $B_\Pi$ is an equitable quotient matrix of $M$, then the eigenvalues of  $B_\Pi$ are also eigenvalues of $M$. Furthermore, if $M$ is nonnegative and irreducible, then $\lambda_{1}(M) = \lambda_{1}(B_\Pi).$
\end{lem}

\begin{lem}(See \cite{D.F})\label{lem::2.2}
Let $n=\sum_{i=1}^t n_i+s$. If $n_{1}\geq n_{2}\geq \cdots\geq n_{t}\geq p$ and $n_{1}<n-s-p(t-1)$, then
$$\rho(K_{s} \nabla (K_{n_{1}}\cup K_{n_{2}}\cup \cdots \cup K_{n_{t}}))<\rho(K_{s} \nabla (K_{n-s-p(t-1)}\cup (t-1)K_{p})).$$
\end{lem}

 Now we shall give a proof of Theorem \ref{thm::1.1}.
\renewcommand\proofname{\bf Proof of Theorem \ref{thm::1.1}}
\begin{proof}

Suppose that $G$ is not a 1-tough graph, there exists some nonempty subset $S$ of $V(G)$ such that $c(G-S)\geq|S|+1$. Let $|S|=s$. Then $G$ is a spanning subgraph of $G_s^1=K_{s} \nabla (K_{n_1}\cup K_{n_2}\cup \cdots \cup K_{n_{s+1}})$ for some integers $n_1\geq n_2\geq \cdots\geq n_{s+1}$ with $\sum_{i=1}^{s+1}n_i=n-s$. Thus,
\begin{equation}\label{equ::1}
\rho(G)\leq\rho(G_s^1),
\end{equation}
where the equality holds if and only if $G\cong G_s^1$.  Then we shall divide the proof into the following three cases.

{\flushleft\bf{Case 1.}} $s\geq\delta+1.$

Let $G_s^2=K_{s} \nabla (K_{n-2s}\cup sK_{1})$. By Lemma \ref{lem::2.2}, we have
\begin{equation}\label{equ::2}
\rho(G_s^1)\leq \rho(G_{s}^2),
\end{equation}
with equality if and only if $(n_1,\ldots,n_{s+1})=(n-2s,1,\ldots,1)$. Observe that $G_s^2$ has the equitable quotient matrix
$$
B_\Pi^s=\begin{bmatrix}
s-1 &s & n-2s\\
s   &0 &0 \\
  s &0 &n-2s-1
\end{bmatrix}.
$$
By a simple calculation, the characteristic polynomial of $B_\Pi^s$ is
\begin{equation}\label{equ::3}
\begin{aligned}
\varphi(B_{\Pi}^s,x)=x^3-(n-s-2)x^2-(s^2+n-s-1)x+s^{2}n-2s^3-s^2.
\end{aligned}
\end{equation}
Also, note that $A(K_{\delta} \nabla (K_{n-2\delta}\cup \delta K_{1}))$ has the equitable quotient matrix $B_{\Pi}^{\delta}$, which is obtained by replacing $s$ with $\delta$ in $B_\Pi^s$. Then
\begin{eqnarray}\label{equ::4}
\varphi(B_{\Pi}^s,x)-\varphi(B_{\Pi}^{\delta},x)=(s-\delta)(x^2+(1-\delta-s)x-2\delta^2
+\delta n-2\delta s+sn-2s^2-\delta-s).
\end{eqnarray}
 Let
$$\tau(x)=x^2+(1-\delta-s)x-2\delta^2+\delta n-2\delta s+sn-2s^2-\delta-s.$$
Then the symmetry axis of $\tau(x)$ is
$$x=\frac{\delta+s-1}{2}<s<n-\delta-1,$$
where the last two inequalities follow from the fact that $s\geq \delta+1$ and $n\geq 2s+1$. This implies that $\tau(x)$ is increasing with respect to $x\geq n-\delta-1$. For $x\geq n\!-\!\delta\!-\!1$, we have

\begin{equation*}
\begin{aligned}
\tau(x)&=x^2\!+\!(1\!-\!\delta\!-\!s)x\!-\!2\delta^2
\!+\!\delta n\!-\!2\delta s\!+\!sn\!-\!2s^2\!-\!\delta\!-\!s\\
&\geq n^2\!-\!(2\delta\!+\!1)n\!+\!(1\!-\!s)\delta\!-\!2s^2~(\mbox{since $x\geq n\!-\!\delta\!-\!1$})\\
&\geq \frac{1}{2}n^2\!-\!\frac{5}{2}\delta n\!+\!\frac{3}{2}\delta\!-\!\frac{1}{2} ~(\mbox{since $ s\leq \frac{n-1}{2}$})\\
&\geq \frac{3\delta}{2}\!-\!\frac{1}{2}~(\mbox{since  ~$n\geq 5\delta$})\\
&> 0~(\mbox{since $\delta\geq 2$}).
\end{aligned}
\end{equation*}
Combining this with (\ref{equ::4}) and $s\geq \delta+1$, we have $\varphi(B_{\Pi}^s,x)>\varphi(B_{\Pi}^{\delta},x)$ for $x\geq n-\delta-1$.
Since $K_{\delta} \nabla (K_{n-2\delta}\cup \delta K_{1})$ contains $K_{n-\delta}$ as a proper subgraph, we have  $\rho(K_{\delta} \nabla (K_{n-2\delta}\cup \delta K_{1}))>\rho(K_{n-\delta})=n-\delta-1$, and so
$\lambda_{1}(B_{\Pi}^s)<\lambda_{1}(B_{\Pi}^{\delta})$. Combining this with Lemma \ref{lem::2.1}, (\ref{equ::1}) and (\ref{equ::2}), we have
$$\rho(G)\leq\rho(G_s^1)\leq\rho(G_s^2)<\rho(K_{\delta} \nabla (K_{n-2\delta}\cup \delta K_{1})).$$

{\flushleft\bf{Case 2.}} $s<\delta$.

Let $G_{s}^{3}=K_{s} \nabla (K_{n-s-(\delta+1-s)s}\cup sK_{\delta+1-s})$. Recall that $G$ is a spanning subgraph of $G_s^1=K_{s} \nabla (K_{n_1}\cup K_{n_2}\cup \cdots \cup K_{n_{s+1}})$, where $n_1\geq n_2\geq \cdots\geq n_{s+1}$ and  $\sum_{i=1}^{s+1}n_i=n-s$. Clearly, $n_{s+1}\geq \delta+1-s$ because the minimum degree of $G_s^1$ is at least $\delta$. By Lemma \ref{lem::2.2}, we have
 \begin{equation}\label{equ::5}
\rho(G_{s}^{1})\leq \rho(G_{s}^{3}),
\end{equation}
where the equality holds if and only if $(n_1,\ldots,n_{s+1})=(n-s-(\delta+1-s)s,\delta+1-s,\ldots,\delta+1-s)$.  Note that $G_s^3$ has the equitable quotient matrix
$$
C_{\Pi'}=\begin{bmatrix}
s-1 &n-s-(\delta-s+1)s & (\delta-s+1)s\\
s   &n-s-(\delta-s+1)s-1 &0 \\
  s &0 &\delta-s
\end{bmatrix}.
$$
By a simple calculation, the characteristic polynomial of $C_{\Pi'}$ is
\begin{equation*}
\begin{aligned}
\phi(C_{\Pi'},x)&=x^{3}\!-\!(\!-\!\delta s\!+\!s^2\!+\!\delta\!+\!n\!-\!2s\!-\!2)x^2\!-\!(\delta^{2}s\!-\!\delta s^2
\!-\!\delta n\!+\!ns\!+\!s^2\!+\!2\delta\!+\!n\!-\!3s\!-\!1)x\\
~~~&~~~\!-\!s(s(\delta\!-\!s)(\delta s\!-\!s^2\!-\!n\!+\!3s)\!+\!\delta^2
\!-\!s\delta\!-\!ns)\!-\!2s^3\!+\!(n\!-\!s)\delta\!-\!ns\!-\!\delta\!+\!s.
\end{aligned}
\end{equation*}
Recall that $\rho(K_{\delta} \nabla (K_{n-2\delta}\cup \delta K_{1}))>n-\delta-1$. Combining these with (\ref{equ::3}), we obtain
\begin{equation*}
\begin{aligned}
&\phi(C_{\Pi'},n\!-\!\delta\!-\!1)-\varphi(B_{\Pi}^\delta,n\!-\!\delta\!-\!1)\\
&=(\delta
\!-\!s)((s\!-\!1)n^2\!+\!(s^2\!-\!(3\delta\!+\!2)s\!+\!3\delta\!+\!1)n
\!+\!
s^4\!-\!(\delta\!+\!3)s^3\!+\!
2s^2\!+\!(2\delta^2\!+\!3\delta)s\!-\!\delta^2\!-\!\delta)\\
&=(\delta\!-\!s)\psi(n),
\end{aligned}
\end{equation*}
where
$$\psi(n)=(s\!-\!1)n^2\!+\!(s^2\!-\!(3\delta\!+\!2)s\!+\!3\delta\!+\!1)n
\!+\!s^4\!-\!(\delta\!+\!3)s^3\!+\!
2s^2\!+\!(2\delta^2\!+\!3\delta)s\!-\!\delta^2\!-\!\delta.$$
Then the symmetry axis of $\psi(n)$ is
$$n=\frac{3\delta-s+1}{2}\leq\frac{3\delta-1}{2}<\frac{2\delta^{2}}{5}+\delta,$$
where the last two inequalities follow from the fact that $s\geq2$ and $\delta\geq3$. This implies that $\psi(n)$ is increasing with respect to $n\geq \frac{2\delta^{2}}{5}+\delta$, and so
\begin{equation*}
\begin{aligned}
\psi(n)&\geq \psi(\frac{2\delta^{2}}{5}+\delta)\\
&=\frac{\delta}{25}[\delta(4(s-1)\delta^{2}-\!10(s-1)\delta\!+\!10s^2\!-\!20s\!+\!35)\!-\!25s^3\!+\!25s^2\!+\!25s]\!+\!s^2(s\!-\!1)(s\!-\!2)\\ &\geq\frac{\delta}{25}[\delta(4s^3\!+\!4s^2\!-\!24s\!+\!41)\!-\!25s^3\!+\!25s^2\!+\!25s]\!+\!s^2(s\!-\!1)(s\!-\!2)
~(~\mbox{since $\delta\geq s+1$})\\
&\geq \frac{\delta}{25}(4s^4-17s^3+5s^2+42s+41)+\!s^2(s\!-\!1)(s\!-\!2)~(~\mbox{since $\delta\geq s+1$})\\
&> 0 ~(~\mbox{since ~$s\geq 2$ and~$\delta\geq 3$}).
\end{aligned}
\end{equation*}
Since $\delta\geq s+1$, it follows that

\begin{equation}\label{equ::6}
\begin{aligned}
\phi(C_{\Pi'},n\!-\!\delta\!-\!1)>\varphi(B_{\Pi}^\delta,n\!-\!\delta\!-\!1).
\end{aligned}
\end{equation}
Next, we take the derivative of $\phi(C_{\Pi'},x)$ and $\varphi(B_{\Pi}^\delta,x)$. Note that $\delta\geq s+1$ and $s\geq 2$. If $x\geq n-\delta-1$, then we obtain that
\begin{equation*}
\begin{aligned}
&\phi'(C_{\Pi'},x)-\varphi'(B_{\Pi}^\delta,x)\\
&=(\delta\!-\!s)((2s\!-\!4)x+
(1\!-\!\delta)s\!+\!n\!+\!\delta\!-\!3) \\
&\geq(\delta\!-\!s)((2s\!-\!4)(n\!-\!\delta\!-\!1)\!+\!(1\!-\!\delta)s\!+\!n\!+\!\delta\!-\!3) ~(\mbox{since $x\geq n\!-\!\delta\!-\!1$, $\delta\geq s+1$ and $s\geq 2$})\\
&=(\delta\!-\!s)((2s\!-\!3)n\!-\!3s\delta+5\delta-s+1)\\
&\geq(\delta\!-\!s)((2s\!-\!3)(\frac{2\delta^{2}}{5}\!+\!\delta)\!-\!3s\delta\!+\!5\delta\!-\!s\!+\!1)~(\mbox{since $n\geq\frac{2\delta^{2}}{5}\!+\!\delta$})\\
&=(\delta\!-\!s)[\frac{\delta}{5}((2s\!-\!3)(2\delta\!-\!3)\!+\!s\!+\!1)\!-\!s\!+\!1]\\
&>0~~(\mbox{since $s\geq 2$ and $\delta\geq s+1$}),
\end{aligned}
\end{equation*}
and so
\begin{equation}\label{equ::7}
\begin{aligned}
\phi'(C_{\Pi'},x)>\varphi'(B_{\Pi}^\delta,x)
\end{aligned}
\end{equation}
for $x\geq n-\delta-1$. We also note that
\begin{equation*}
\begin{aligned}
&\phi'(C_{\Pi'},n\!-\!\delta\!-\!1)\\
&=n^2\!+\!(2\delta s\!-\!2s^2\!-\!5\delta\!+\!3s\!-\!1)n\!-\!3\delta^{2} s\!+\!3\delta s^2\!+\!5\delta^2\!-\!6s\delta\!+\!s^2\!+\!2\delta\!-\!s\\
&\geq\frac{\delta}{25}[\delta(4\delta^2\!+\!(20s\!-\!30)\delta\!-\!20s^2\!+\!5s\!+\!15)\!+\!25s^2\!-\!75s\!+\!25]\!+\!s^2\!-\!s ~(\mbox{since $n\geq\frac{2\delta^{2}}{5}\!+\!\delta$})\\
&\geq\frac{\delta}{25}(4s^3\!+\!32s^2\!-\!83s\!+\!14)\!+\!s^2\!-\!s ~(\mbox{since $\delta\geq s+1$})\\
&> 0 ~(\mbox{since $\delta\geq 3$ and ~$s\geq 2$}),
\end{aligned}
\end{equation*}
and the symmetry axis of $\phi'(C_{\Pi'},x)$ is
$x=\frac{-\delta s+s^2+\delta+n-2s-2}{3}<n-\delta-1$ due to $\delta\geq s+1$ and $n\geq \frac{2}{5}\delta^{2}+\delta$. It follows that $\phi(C_{\Pi'},x)$ is increasing with respect to $x\geq n-\delta-1$. Combining this with (\ref{equ::6}) and (\ref{equ::7}), we can deduce that $\lambda_{1}(C_{\Pi'})<\lambda_{1}(B_{\Pi}^{\delta})$. Furthermore, according to Lemma \ref{lem::2.1}, (\ref{equ::1}) and (\ref{equ::5}), we have
$$\rho(G)\leq\rho(G_s^1)\leq\rho(G_s^3)<\rho(K_{\delta} \nabla (K_{n-2\delta}\cup \delta K_{1})).$$

{\flushleft\bf{Case 3.}} $s=\delta.$

By Lemma \ref{lem::2.2}, we have
$$\rho(G_{s}^1)\leq \rho(K_{\delta} \nabla (K_{n-2\delta}\cup \delta K_{1})),$$
with equality holding if and only if $G_{s}^1\cong K_{\delta} \nabla (K_{n-2\delta}\cup \delta K_{1})$. Combining this with (\ref{equ::1}), we may conclude that
$$\rho(G)\leq\rho(K_{\delta} \nabla (K_{n-2\delta}\cup \delta K_{1})),$$
where the equality holds if and only if $G\cong K_{\delta} \nabla (K_{n-2\delta}\cup \delta K_{1})$.

 This completes the proof. \end{proof}

\section{Proof of Theorem \ref{thm::1.2}}

In this section, we give the proof of Theorem \ref{thm::1.2}. The following lemmas are used in the sequel.
\begin{lem}(See \cite{C.H})\label{lem::3.1}
Let $G$ be a 1-tough graph with the degree sequence $(d_{1}, d_{2},\ldots, d_{n})$, where $d_{1}\leq d_{2}\leq\cdots\leq d_{n}$. If for all $1\leq i<\frac{n}{2}$, $d_{i}\leq i$ implies $d_{n-i+1}\geq n-i$,
then $G$ contains a hamiltonian cycle.
\end{lem}

\begin{lem}(See \cite{HSF,V.N})\label{lem::3.2}
Let $G$ be a graph on $n$ vertices and $m$ edges with $\delta\geq 1$. Then
$$\rho(G) \leq \frac{\delta-1}{2}+\sqrt{2 m-n \delta+\frac{(\delta+1)^{2}}{4}},$$
with equality if and only if $G$ is either a $\delta$-regular graph or a bidegreed graph
in which each vertex is of degree either $\delta$ or $n-1$.
\end{lem}

\begin{lem}(See \cite{V.N})\label{lem::3.3}
For nonnegative integers $p$ and $q$ with $2q \leq p(p-1)$ and $0 \leq x \leq p-1$, the function $f(x)=(x-1) / 2+\sqrt{2 q-p x+(1+x)^{2} / 4}$ is decreasing with respect to $x$.
\end{lem}

\begin{lem}(See \cite{A.A})\label{lem::3.4}
Let $G$ be a connected 1-tough graph of order $n\geq3$ with the degree sequence $(d_{1}, d_{2},\ldots, d_{n})$, where $d_{1}\leq d_{2}\leq\cdots\leq d_{n}$, and let $\theta$ be the smallest integer for which $d_{\theta}\geq \frac{n}{2}$. If $1 \leq \theta \leq \delta+2$, then either $G$ contains a hamiltonian cycle or $G\cong M_{n}$.
\end{lem}

For any $S\subseteq V(G)$, let $G[S]$ be the subgraph of $G$ induced by $S$, and let $e(S)$ be the number of edges in $G[S]$. Denote by $G-v$ and $G-uv$ the graphs obtained from $G$ by deleting the vertex $v\in V(G)$ and the edge $uv\in E(G)$, respectively. Similarly, $G+uv$ is obtained from $G$ by adding the edge $uv\notin E(G)$.

\begin{lem}\label{lem::3.5}
Suppose that $G$ is a graph on $n\geq k+2$ vertices. If $e(G)\geq{n-1\choose 2}+k+2$, then for edge set $M$ of size $k$ with $\Delta(G[V(M)])\leq 2$, $G$ has a hamiltonian cycle containing $M$.
\end{lem}

 \noindent \textbf{Proof.}
Suppose that the result does not hold. Let $G$ be a maximal counterexample. Then $G$ has a hamiltonian  path $P$ containing $M$, say $P=x_1x_2\ldots x_n$. By hypothesis, we have $x_1x_n\notin E(G)$. Let $S$ be a subset of $V(M)$ such that $x\in S$ is the second vertex incident with $e$ along $P$ for each $e\in M$.
Let $I_1=\{j-1\ |\ x_j\in N_G(x_1)\}$,  $I_2=\{j\ |\ x_j\in N_G(x_n)\}$ and $I_3=\{j-1\ |\ x_j\in S\}$. Then $n\notin I_1\cup I_2$.
Since $e(G)\geq {n-1\choose 2}+k+2$, it follows that $d_G(x_1)+d_G(x_n)\geq n+k$. This implies that $|I_1|+|I_2|\geq  n+k$ and $|I_1\cap I_2|\geq k+1.$  Note that  $|I_3|=|S|=k$. Thus there exists some $j\in (I_1\cap I_2)\backslash I_3$ such that $x_{j+1}\in N_G(x_1)$, $x_{j}\in N_{G}(x_n)$ and $x_jx_{j+1}\notin M$, and so we obtain a hamiltonian cycle $C$ containing $M$ of $G$ from $(E(P)-x_jx_{j+1})\cup \{x_jx_n,x_{j+1}x_1\}$, a contradiction. This completes the proof. \qed

\begin{lem}(See \cite{QZ})\label{lem::3.6}
Let $G$ be a connected graph of order $n\geq 13$ with $\delta(G)\geq3$. If
$$e(G)\geq{n-2\choose 2}+6,$$
then $G$ is hamiltonian-connected, unless $G\cong K_{3}\nabla(K_{n-5}\cup2K_{2})$.
\end{lem}

Let $E_1=\{u_{1}v_{1},u_{2}v_{2},\ldots,u_{p}v_{p}\}$ where $p\geq 1$, and let $G(E_1)
$ be the graph with vertex set $V(G)$ and edge set $E_1\cup E(G)$. Particularly, $G(E_1)=G$ if $E_1\subseteq E(G)$.

\begin{lem}\label{lem::3.7}
Suppose that $G$ is a connected graph on $n\geq 15$ vertices. If $e(G)\geq{n-1\choose 2}$ and $\delta(G)\geq 4$, then for any $e,f\in E(G)$, $G$ has a hamiltonian cycle containing $e,f$.
\end{lem}
\noindent \textbf{Proof.}
Let $e,f\in E(G)$, say $e=v_1v_2$ and $f=u_1u_2$, and let $S:=\{v_1,v_2,u_1,u_2\}$.
Since $e(G)\geq {n-1\choose 2}>{n-2\choose 2}+6$ for $n\geq 15$, by Lemma \ref{lem::3.6}, $G$ is hamiltonian-connected.
 So $G$ has a hamiltonian cycle $C$ containing $e$. If $C$ contains $f$, then $C$ is a desired cycle. Now we assume that $f\notin E(C)$. Thus we may obtain a hamiltonian path $P$ containing $e$ and $f$ from $E(C)\cup \{f\}$.  Without loss of generality, we assume that $P=x_1x_2\ldots x_n$ is a hamiltonian path of $G$ from $x_1$ to $x_n$, where $x_1x_n\notin E(G)$ and $d_G(x_1)\leq d_G(x_n)$. Suppose that $u_2$ and $v_2$ appear after $u_1$ and $v_1$ along $P$, respectively. Let $S':=\{u_2,v_2\}$, and let $i=|N_G(x_1)\cap S'|$. Then we divide the proof into the following two cases.

{\flushleft\bf{Case 1.}} $d_G(x_1)+d_G(x_n)\geq n+i$.

Let $I_1=\{j-1\ |\ x_j\in N_G(x_1)\}$ and let $I_2=\{j\ |\ x_j\in N_G(x_n)\}$. Then $n\notin I_1\cup I_2$. Since $|I_1|+|I_2|\geq n+i$, it follows that $|I_1\cap I_2|\geq i+1$. This implies there exists some $j\in I_1\cap I_2$ such that $x_j\in N_G(x_n)$, $x_{j+1}\in N_G(x_1)$ and $x_jx_{j+1}\notin \{e,f\}$. Then we obtain a hamiltonian cycle containing $e$ and $f$ of $G$ from $(E(P)-x_{j}x_{j+1})\cup \{x_1x_{j+1},x_{n}x_j\}$.

\medskip
{\flushleft\bf{Case 2.}} $d_G(x_1)+d_G(x_n)\leq n+i-1$.

Note that $i\leq 2$ and  $d_G(x_1)\leq d_{G}(x_n)$. Then $d_G(x_1)\leq (n+i-1)/2\leq (n+1)/2$, and so
\begin{align}\label{equ::8}
e(G-x_1)\geq {n-1\choose 2}-(n+1)/2={n-2\choose 2}+(n-5)/2\geq {n-2\choose 2}+ 5
\end{align}
by $n\geq 15$. Note that $x_1$ is an end vertex of $P$. Then $x_1\notin \{v_2,u_2\}$. Without loss of generality, we may assume that $x_1\neq u_1$.
 We consider the following two subcases.

{\flushleft\bf{Subcase 2.1.}} $e$ and $f$ are not incident.

If $x_1=v_1$, since $d_{G}(x_1)\geq 4$, there exists $y\in N_G(x_1)\backslash\{v_2\}$. Combining this with Lemma \ref{lem::3.5} and (\ref{equ::8}), $G(v_2y)-x_1$ has a hamiltonian cycle $C$ containing $f$ and $v_2y$. Then $(C-v_2y)+e+x_1y$ is a hamiltonian cycle of $G$ containing $e$ and $f$. Now we may assume that $x_1\neq v_1$. Since $d_G(x_1)\geq 4$, there exist $y_1,y_2\in N_G(x_1)$ such that $y_1y_2\notin \{e,f\}$. Again by Lemma \ref{lem::3.5}, we deduce that $G(y_1y_2)-x_1$ has a hamiltonian cycle $C$ containing $e,f$ and $y_1y_2$. By adding edges $y_2x_1,y_1x_1$ to $C-y_1y_2$, we obtain a hamiltonian cycle of $G$ containing $e$ and $f$.

{\flushleft\bf{Subcase 2.2.}} $e$ and $f$ are incident (i.e., $v_2=u_1$).

If $x_1=v_1$, then $e=x_1x_2$. Let $y\in N_G(x_1)\backslash S$. Then $G(x_2y)-x_1$ has hamiltonian cycle $C$ containing $x_2y$ and $f$ due to Lemma \ref{lem::3.5} and (\ref{equ::8}), and so we obtain a hamiltonian cycle of $G$ containing $e,f$ from $(C-x_2y)+e+x_1y$. Next we assume that $x_1\neq v_1$, i.e., $x_1\notin S$.
Since $d_{G}(x_1)\geq 4$ and $N_G(x_1)\neq S$, there exist $y_1\in N_G(x_1)\backslash S$ and $y_2\in N_G(x_1)\backslash\{x_2,y_1\}$. Again by Lemma \ref{lem::3.5}, we deduce that $G(y_1y_2)-x_1$ has a hamiltonian cycle $C$ containing $e,f$ and $y_1y_2$.
Then we may obtain a hamiltonian cycle of $G$ containing $e$ and $f$ from $(C-y_1y_2)+x_1y_1+x_1y_2$.

This completes the proof. \qed

\renewcommand\proofname{\bf Proof of Theorem \ref{thm::1.2}}
\begin{proof}
Suppose that a 1-tough graph $G$ contains no hamiltonian cycles and its degree sequence is $(d_{1}, d_{2},\ldots, d_{n})$, where $\delta=d_{1}\leq d_{2}\leq\cdots\leq d_{n}$. By Lemma \ref{lem::3.1}, there is an integer $k < \frac{n}{2}$ such that $d_{k} \leq k$ and $d_{n-k+1}\leq n-k-1$. Since $\delta\geq 2$ and $k\geq d_{k}\geq \delta$, it follows that

\begin{equation}\label{equ::9}
\begin{aligned}
2m&=\sum_{i=1}^{n}d_{i}\\
 &\leq k^{2}+(n-2k+1)(n-k-1)+(k-1)(n-1)\\
 &=n^{2}-(2k+1)n+3k^{2}.
\end{aligned}
\end{equation}
By Lemmas \ref{lem::3.2}, \ref{lem::3.3}, and the fact $\delta\geq 2$, we obtain

\begin{equation}\label{equ::10}
\rho(G)\leq \frac{1}{2}+\sqrt{2m-2n+\frac{9}{4}}.
\end{equation}
Note that $\rho(G)\geq\rho(M_{n})>\rho(K_{n-3})=n-4$. Combining this with (\ref{equ::10}), we have
\begin{equation}\label{equ::11}
2e(G)=2m > n^{2}-7n+18,
\end{equation}
and so $e(G)\geq{n-3\choose 2}+4$. Recall that $k< \frac{n}{2}$. Then $n\geq 2k+1$. For $k\geq 4$, by (\ref{equ::9}) and (\ref{equ::11}), we deduce that
\begin{equation*}
(2k-6)(2k+1)\leq (2k-6)n< 3k^{2}-18,
\end{equation*}
from which we have $k<\frac{10+\sqrt{52}}{2}$, and so $ k\leq 8$. Thus

$$n<\frac{3k^2-18}{2k-6}=\frac{1}{2}(3k+9+\frac{9}{k-3})<18$$
for $4\leq k\leq 8$, which is impossible because $n\geq 18$.
We consider $k=2$ or $k=3$ in the following. If $G$ is a subgraph of $M_{n}$, then
$\rho(G)\leq\rho(M_{n})$, where the equality holds if and only if $G\cong M_{n}$,
which is impossible because $\rho(G)\geq\rho(M_{n})$ and $G\ncong M_{n}$. Assume that $G$ is not a subgraph of $M_{n}$. Let $\theta\geq 1$ be the smallest integer for which $d_{\theta}\geq \frac{n}{2}$ in $G$. We divide the proof into the following two cases.

{\flushleft\bf{Case 1.}} $k=2.$

Note that $\delta\geq 2$ and $d_{2}\leq 2$. Then $d_{1}=d_{2}=\delta=2$. Let $v_{1}, v_{2}$ and $v_{3}$ be the vertices of $G$ such that $d_{G}(v_{i})=d_i$ for $1\leq i\leq 3$. If $d_{3}\leq 5$, since $e(G)\geq{n-3\choose 2}+4$, we get
 \begin{equation*}
e(V(G)-\{v_{1},v_{2},v_{3}\})\geq \binom{n-3}{2}+4-9=\binom{n-3}{2}-5,
\end{equation*}
and so
\begin{equation*}
d_{4}\geq\binom{n-3}{2}-5-\binom{n-4}{2}=n-9\geq\frac{n}{2}
\end{equation*}
for $n\geq 18$, which implies that $\theta=4$. Note that $\delta=2$, and so $\theta-\delta= 2$. Then $G$ contains a hamiltonian cycle by Lemma \ref{lem::3.4}, a contradiction. Thus we assume that $d_{3}\geq 6$, and so  $\delta(G\!-\!v_1\!-\!v_2)\geq 4$. Since $e(G)\geq{n-3\choose 2}+4$, it follows that
\begin{equation}\label{equ::12}
e(V(G)-\{v_{1},v_{2}\})\geq {n-3\choose 2}.
\end{equation}
 If $v_{1}v_{2}\in E(G)$, then $N_{G}(v_1)\cap N_{G}(v_2)=\emptyset$. Otherwise, we assume that $N_{G}(v_1)\cap N_{G}(v_2)=\{u\}$. Thus $u$ is a cut vertex, which is impossible because $G$ is a 1-tough graph. Let $w_{1}\in N_{G}(v_1)\backslash\{v_{2}\}$ and $w_{2}\in N_{G}(v_2)\backslash\{v_{1}\}$. Combining (\ref{equ::12}) with Lemma \ref{lem::3.7}, $G(w_{1}w_{2})-v_1-v_2$ has a hamiltonian cycle $C_1$ containing $w_{1}w_{2}$. Thus we can deduce that $G$ has a hamiltonian cycle from $C_{1}-w_{1}w_{2}+w_{1}v_{1}+v_{1}v_{2}+v_{2}w_{2}$, a contradiction. If $v_{1}v_{2}\notin E(G)$, we assume that $N_{G}(v_1)=\{w_{1},w_{2}\}$ and $N_{G}(v_2)=\{w_{3},w_{4}\}$. By using the similar analysis as above, $G(w_{1}w_{2},w_{3}w_{4})-v_1-v_2$ has a hamiltonian cycle $C_{2}$ containing $w_{1}w_{2}$ and $w_{3}w_{4}$, and so  $C_{2}-w_{1}w_{2}-w_{3}w_{4}+w_{1}v_{1}+w_{2}v_{1}+w_{3}v_{2}+w_{4}v_{2}$ is a hamiltonian cycle of $G$, which also leads to a contradiction.

{\flushleft\bf{Case 2.}} $k=3.$

Since $d_{3}\leq 3$, there exist three vertices $v_{1}, v_{2}$ and $v_{3}$ such that $d_{G}(v_{i})\leq 3$ for $1\leq i\leq3$. Combining this with (\ref{equ::11}), we get
 \begin{equation*}
e(V(G)-\{v_{1},v_{2},v_{3}\})\geq \binom{n-3}{2}+4-9=\binom{n-3}{2}-5,
\end{equation*}
and so
\begin{equation*}
d_{4}\geq\binom{n-3}{2}-5-\binom{n-4}{2}=n-9\geq\frac{n}{2}
\end{equation*}
for $n\geq 18$, which implies that $\theta=4$. Note that $2\leq\delta\leq 3$, and so $\theta-\delta\leq 2$. Therefore, $G$ contains a hamiltonian cycle by Lemma \ref{lem::3.4}, a contradiction.

This completes the proof.\end{proof}

\section{Proof of Theorem \ref{thm::1.3}}

\begin{lem}(See \cite{Y.Hong})\label{lem::4.1}
Let $G$ be a graph with $n$ vertices and $m$ edges. Then
                  $$\rho(G)\leq\sqrt{2m-n+1},$$
where the equality holds if and only if $G$ is a star or a complete graph.
\end{lem}

 Now we shall give a short proof of Theorem \ref{thm::1.3}.
\renewcommand\proofname{\bf Proof of Theorem \ref{thm::1.3}}
\begin{proof}

Suppose that $G$ is not a $t$-tough graph, there exists some nonempty subset $S$ of $V(G)$ such that $tc(G-S)>|S|$. Let $|S|=s$ and $c(G-S)=c$. Then $G$ is a spanning subgraph of $G_1=K_{tc-1} \nabla (K_{n_1}\cup K_{n_2}\cup \cdots \cup K_{n_{c}})$ for some integers $n_1\geq n_2\geq \cdots\geq n_{c}$ with $\sum_{i=1}^{c}n_i=n-tc+1$. Thus,
\begin{equation}\label{equ::13}
\rho(G)\leq\rho(G_1),
\end{equation}
where the equality holds if and only if $G\cong G_1$.
Let $G_2=K_{tc-1} \nabla (K_{n-(t+1)c+2}\cup (c-1)K_{1})$. By Lemma \ref{lem::2.2}, we have
\begin{equation}\label{equ::14}
\rho(G_1)\leq \rho(G_2),
\end{equation}
with equality if and only if $(n_1,\ldots,n_{c})=(n-(t+1)c+2,1,\ldots,1)$.

If $c=2$, then $G_2\cong K_{2t-1} \nabla (K_{n-2t}\cup K_{1})$. Combining this with (\ref{equ::13}) and (\ref{equ::14}), we may conclude that
$$\rho(G)\leq\rho(K_{2t-1} \nabla (K_{n-2t}\cup K_{1})),$$
where the equality holds if and only if $G\cong K_{2t-1} \nabla (K_{n-2t}\cup K_{1})$.
For $c\geq3$, by Lemma \ref{lem::4.1}, we have
\begin{equation}\label{equ::15}
\begin{aligned}
\rho(G_2)&\leq\sqrt{2e(G_2)-n+1}\\
&=\sqrt{[(n-c+1)(n-c)+2(tc-1)(c-1)]-n+1}\\
&=\sqrt{(1+2t)c^{2}-(2n+2t+3)c+n^{2}+3}.
\end{aligned}
\end{equation}
Let $f(c)=(1+2t)c^{2}-(2n+2t+3)c+n^{2}+3$. Since $n\geq (t+1)c-1$, we have $3 \leq c\leq(n+1)/(t+1)$. By a simple calculation,
\begin{equation*}
\begin{aligned}
f(3)-f\Big(\frac{n+1}{t+1}\Big)=\frac{(n-3t-2)(n-4t^2-6t-1)}{(t+1)^{2}}>0,
\end{aligned}
\end{equation*}
where the inequality follows from the fact that $n\geq 4t^2+6t+2$. This implies that, for $3\leq c\leq(n+1)/(t+1)$, the maximum value of $f(c)$ is attained at $c=3$, and so from (\ref{equ::15}), we deduce that
\begin{equation}\label{equ::16}
\begin{aligned}
\rho(G_2)&\leq \sqrt{f(3)}\\
&=\sqrt{(n-2)^{2}-(2n-12t+1)}\\
&\leq \sqrt{(n-2)^{2}-(2(4t^2+6t+2) - 12t+1)} ~(\mbox{since $n\geq 4t^2+6t+2$})\\
&< n-2~(\mbox{since $t\geq 1$}).\\
\end{aligned}
\end{equation}
Since $K_{2t-1} \nabla (K_{n-2t}\cup K_{1})$ contains  $K_{n-1}$ as a proper subgraph, we have
 $$\rho(K_{2t-1} \nabla (K_{n-2t}\cup K_{1}))>\rho(K_{n-1})=n-2.$$
Combining this with (\ref{equ::13}), (\ref{equ::14}) and (\ref{equ::16}), we have $$\rho(G)\leq \rho(G_1)\leq \rho(G_2)<\rho(K_{2t-1} \nabla (K_{n-2t}\cup K_{1})).$$

This completes the proof. \end{proof}

\begin{remark}
For two positive integers $t$ and $c$, one can verify that $K_{tc-1}\nabla (K_{n-(t+1)c+2}\cup (c-1)K_{1})$ is not a $t$-tough graph. Notice that $K_{2t-1}\nabla(K_{n-2t}\cup K_{1})$ may not be the extremal graph in Theorem \ref{thm::1.3} if the value of $n-t$ is not large enough. For example, by Matlab programming, we have $\rho(K_{17}\nabla 3K_1)=18.72381$ and $\rho(K_{11}\nabla(K_{8}\cup K_{1}))=18.35161$, and so $\rho(K_{17}\nabla 3K_1)>\rho(K_{11}\nabla(K_{8}\cup K_{1}))$.
Thus, we consider the order of $G$ is sufficiently large with respect to $t$ in Theorem \ref{thm::1.3}.
\end{remark}

\end{document}